\documentclass[12pt]{amsart}
\usepackage{graphicx}
\usepackage{amsmath}

\newcommand {\PP}{{I\kern-.3em P}}
\newcommand {\ZZ}{{Z\kern-.45em Z}}
\newcommand {\RR}{{I\kern-.3em R}}

\newcommand {\HH}{{\mathbb{H}}}
\newcommand {\NN}{{I\kern-.3em N}}

\newcommand{\beq}{\begin{equation}}
\newcommand{\eeq}{\end{equation}}

\def\d{\displaystyle}

\def\l{$\lambda\,\,$}
\def\s{$S\,\,$}
\def\p{$p\,\,$}
\def\pt{$\tilde{p}\,\,$}
\def\gam{$\gamma\,\,$}
\def\gamt{$\tilde{\gamma}\,\,$}
\newtheorem{thm}{Theorem}
\newtheorem{lemma}{Lemma}
\newtheorem{prop}{Proposition}
\newtheorem{cor}{Corollary}
\theoremstyle{remark}

\theoremstyle{property}

\theoremstyle{definition}

\begin{document}
\title{The distribution of geodesic excursions into the neighborhood of a cone singularity on a hyperbolic 2-orbifold }
\author{Andrew Haas }

\address{Department of Mathematics,
The University of Connecticut,
Storrs, CT. 06269-3009}
\begin{abstract}
A generic geodesic on a finite area, hyperbolic 2-orbifold exhibits an infinite sequence of penetrations into a neighborhood of a cone singularity of order $k\geq 3,$ so that the sequence of depths of maximal penetration has a limiting distribution. The distribution function is the same for all such surfaces and is described by a fairly simple formula.

 \end{abstract}

\email{haas@math.uconn.edu}
\date{}
\maketitle
\markboth{The distribution of geodesic excursions}
{The distribution of geodesic excursions  }

\section{introduction}
Many of the concepts of number theory can be realized in a geometric setting. This is especially true in the field of classical diophantine approximation, where much of the theory has been reformulated   as statements about the depth of penetration of a geodesic into the cusp, or non-compact, end of   the Modular Surface.

The Theorem of Bosma, Jager and Weidijk \cite{bjw}, called the Lenstra Conjecture, describes the distribution of  values $\theta_n=q_n|q_nx-p_n|$, for almost all real numbers $x$, where $p_n/q_n$ is the $n^{th}$ continued fraction convergent to $x$. Bosma proved a closely related result for approximation by mediants \cite{bosma}.
In \cite{haas} we proved analogues of these results which describe the distribution of the depth of maximal penetration  of a generic geodesic into the cusp end of a finite area hyperbolic 2-orbifold. These purely geometric results were then applied to prove   Fuchsian group versions of \cite{bjw} and \cite{bosma}.  The distribution of the $\theta_n,$ defined with respect to a finite area Fuchsian group, is given independent of the group and agrees with the classical result.

In this paper we look at how a generic geodesic on an orbifold surface "approximates" a cone singularity of the surface. This type of investigation was begun in \cite{dukehaas}. In that paper, the geometric Markoff theory was generalized to this setting. Here it is shown that the sequence of depths of  maximal penetration of a generic geodesic into a neighborhood of a cone singularity has a limiting density which is explicitly computable and the same for all finite area surfaces and most orbifold surfaces. As in \cite{haas}, the resulting densities  bear a striking similarity to the ones found in \cite{bjw} and \cite{bosma}, although they are somewhat more complex and arcane. The similarity becomes particularly clear when the distributions are expressed in terms of area.

\subsection{The main results}\label{section1.1}
We shall begin with a heuristic treatment of some of the definitions and then state the main results of the paper.
 A hyperbolic 2-orbifold  $S$ is represented as the quotient of the hyperbolic plane $\HH$ by a Fuchsian group $G.$ We shall assume that \s has finite area and, if \s is compact, $G$ is not a triangle group.
A point  $p$ on \s is a cone point  of order $k\geq 3$ if,    under the projection from   $\HH$ to \s,  $p$ is covered by a point \pt which has non-trivial stabilizer of order $k$ in $G$.  

An $r$-excursion $e,$ of a geodesic ray \gam into a neighborhood of the cone point,
is an arc of \gam associated to a point of maximal penetration of the ray \gam into the radius $r$ neighborhood, $B_r(p)$, of the cone point. 
The depth of the excursion $e$, denoted $d(e)$, is the "distance" from $e$ to \p.
 When $r$ is sufficiently small, the excursions are simply the arcs $\gamma \cap B_r(p)$ and $d(e)$ is the actual distance. The precise definition becomes more complicated with larger values of $r$. An $r$-excursion $e,$ which is associated with a closed arc of \gam that loops about $p$, is called an approximating $r$-excursion. $e$ is called an approximating excursion if it is an  approximating $r$-excursion for some value $r$. These definitions parallel  the geometric analogues of mediant and continued fraction approximations, as they appear on the classical Modular surface.

 The $r$-excursions of a ray \gam naturally form a sequence $\{e_j\}$, which goes out the end of \gam if the sequence is infinite. In a similar fashion the approximating excursions of \gam form a sequence $\{e^{*}_j\}$. For a generic ray $\gamma$,   both   sequences are infinite for all $r>0$.
 
Let $\# A$ denote the cardinality of the set $A$. For  $r>0,$ $0\leq z \leq r$  and a  geodesic ray \gam with infinitely many  $r$-excursions $e_j$ into the radius $r$ neighborhood of a cone point \p of order $k$, the distribution of the depths of $r$-excursions is defined as
  \begin{equation}\label{Edist}
  dist_k(r,z)(\gamma)= \lim_{n\rightarrow\infty} \frac{1}{n}\#\{j\,|\,1\leq j\leq n , d(e_j)\leq z\}, 
  \end{equation}
  if the limit exists.
 Similarly, if  $e_j^{*}$ is an infinite sequence of approximating excursions of \gam into the neighborhood of a cone point \p of order $k$, then the distribution,  $dist_k^{*}(z)(\gamma),$  of  the depths of approximating excursions is defined just as the function (\ref{Edist}) above, with $e_j^*$ replacing $e_j$ in the definition, and  $0\leq z\leq r_k$, where $r_k=  \sinh^{-1}(\cot\frac{\pi}{k}        ). $ Note that the function does not depend on a parameter $r$.
 In Corollary \ref{cor2}, of Section \ref{section42}, we show that every approximating excursion is an $r$-excursion for some  
  $r< r_k$.
 

  We are now in a position to state the main results of the paper.  
\begin{thm}\label{thm1}
For almost all geodesic rays  \gam, for $r>0$ and   $0\leq z\leq r$ the distribution $dist_k(r,z)(\gamma)$ converges to a function $dist_k(r,z).$ When $r\leq   r_k ,$ 
$$ dist_k(r,z) =\d{\frac{  \sinh z }{\sinh r}} .$$
  
\noindent When $r >  r_k ,$
    
 $$   dist_k(r,z) =
  \left\{\begin{array}{cc}
  \d  { \frac { \frac{\pi}{k}\sinh z }{\sinh r \tan^{-1}(\frac{1}{\sinh r}) + \log(\sin  \frac{\pi}{k}\cosh r)}  }
       & {\rm if} \,  z\leq r_k\\  
    \\
   \d{ \frac{\sinh z \tan^{-1}(\frac{1}{\sinh z}) + \log(\sin  \frac{\pi}{k}\cosh z)}
   {\sinh r \tan^{-1}(\frac{1}{\sinh r}) + \log(\sin  \frac{\pi}{k}\cosh r)}}
     & {\rm if} \,  z>r_k
       \end{array} \right.$$
   \end{thm}

 The result for approximating excursions is given in
 
 \begin{thm}\label{thm2}
For almost all geodesic rays  \gam  and for $0\leq z\leq r_k$ the distribution $dist_k^{*}(z)(\gamma)$ converges to a function $dist_k^{*}(z)$, which  is given by the
formula
 \begin{equation*}\label{formula2} 
  dist_k^{*}(z)= 
\left\{\begin{array}{cc}
\d{\frac{\frac{\pi}{k}\sinh z}{\log (2\cos\frac{\pi}{k} )}}
 & {\rm if} \, z\leq\delta_k\\
 \\
\d{ \frac{\sinh z\tan^{-1}(\varphi_k(\sinh z)) +\log(2\cosh z \sin \frac{\pi}{k} \cos\frac{\pi}{k})}{\log (2\cos\frac{\pi}{k}) } }
  & {\rm if} \, z>\delta_k
         \end{array} \right .
  \end{equation*}
  where   $\delta_k=\sinh^{-1}(\cot  \frac{2\pi}{k}  ) $
  and $\varphi_k(x)= \frac{1-x \tan  \frac{\pi}{k}}{x+ \tan  \frac{\pi}{k}}.$  When $k=3\, \text{and}\, 4,\, \delta_k\leq 0$ and the first case does not occur.
   \end{thm}
  
 \subsection{Outline of the paper}
 
 There are two main pieces to the proofs of Theorems \ref{thm1} and \ref{thm2}. First in Theorem \ref{reduction}, we prove the existence of the distributions and describe their values by certain integrals. In Theorem \ref{thmlambda}, we complete the proof by computing the values of these integrals. 
 
 Section \ref{section2} is concerned with background material and a precise treatment of excursions.   The statement   of  Theorem \ref{reduction} is   in Section \ref{section3} and the proof occupies the remainder of Section \ref{ergodic}.  Theorem \ref{thmlambda} is the focus of Section  \ref{section40}. In the last section, Section \ref{section5}, we see how the results look if area, rather than distance, is used to define the depth of an excursion. Written in terms of area, it is easier to see how the work with cone points relates to the more classical case.

\section{ Excursions into a cone neighborhood  }\label{section2}

\subsection{Basic notions and definitions  }
The quotient of the Poincar\'{e} upper-half plane $\HH$ by a Fuchsian group $G$ is called a hyperbolic 
2-orbifold. Let $\pi:\HH\rightarrow \HH/G=S$ denote the projection from the upper-half plane to the orbifold $S$.
A point $p\in S$ is called a {\it cone point of order k} if there is a primage \pt of \p in $\HH$ so that the stabilizer of \pt in $G$ is generated by a transformation $T_k$ of 
order $k$.  

Throughout  this paper we shall assume that $S$ is a finite area hyperbolic 2-orbifold with a cone point \p of order $k\geq 3$ and that, if \s is compact then $G$ is not a triangle group.
One important consequence of this last hypothesis is that,  if \s is compact, then there is a simple closed geodesic on \s that does not pass through \p. This follows from the fact that such a surface will contain a homotopically non-trivial, simple loop in the complement of the cone points  on \s, which does not go around a cone point
We should mention that in the special case where $G$ is a group with signature $(0;2,2,2,k),$ the simple closed geodesic degenerates to a geodesic ray going back and forth between order two cone points.

For $r> 0,$ let $\tilde{B}_r(\tilde{q})$ denote the open hyperbolic disc in $\HH$ of radius $r$ centered at $\tilde{q}.$ When $r=0$   set $\tilde{B}_0(\tilde{q})=\{\tilde{q}\}$. For \pt covering the cone point \p, 
$\tilde{B}_r(\tilde{p})$ projects to what we shall call the {\it cone neighborhood of p of radius r}, written $B_r(p)$. When $G$ is not a triangle group, if $ r <r_k$ then  for $g\in G$,
$g(\tilde{B}_r(\tilde{p}))\cap \tilde{B}_r(\tilde{p})\not =\emptyset$ if and only if $g$ is in the stabilizer of \pt, \cite{Beardon}.   Consequently, the 
projection of $\tilde{B}_r(\tilde{p})$ is exactly $k$-to-1 in the complement of $\tilde{p}$.  For larger values of $r$, the projection of $\tilde{B}_r(\tilde{p})$ can be considerably more complicated.

Let   $\gamma:(-\infty,\infty)\rightarrow S$  be a geodesic  parameterized by arc length. A lift $\tilde{\gamma}$ of \gam to $\HH$
  has endpoints $\tilde{\gamma}_+$ and   $\tilde{\gamma}_-$
in the extended real line $\widehat{\RR}=\RR \cup \{\infty \}$, representing respectively the limits at infinity and
minus infinity. If the domain of $\gamma$ is restricted to $[0,\infty)$ then we shall refer to it as a geodesic ray. A lift 
of the ray $\tilde{\gamma}$ has the single endpoint $\tilde{\gamma}_+$ at infinity. Henceforth, all geodesics shall   be parameterized by arc length.

\subsubsection{The distinguished ray $\lambda$}
It is necessary to choose a distinguished, simple geodesic ray $\lambda$ on \s with initial point $\lambda(0)=p,$ which is used to 
catalogue the excursions of a geodesic on \s into a cone neighborhood of \p. The two cases, \s non-compact and \s
compact, are considered separately. 

If \s is non-compact then it has a cusp end, also called a puncture, and there is at least one simple geodesic
rays with initial point \p that eventually go out the cusp end. Choose  such  a ray and call it $\lambda$. 

Now suppose  \s is compact.   We have assumed that $G$ is not a triangle group and therefore, on \s there is a simple closed geodesic $\beta$ that does not contain \p.
We shall choose $\lambda$ to be a simple geodesic ray that twirls asymptotically into $\beta.$ It can be constructed as follows. Let
$\delta$ be a minimal length, and therefore simple, geodesic arc with its initial  point at \p and its terminal point on $\beta.$
Let $\tilde{\beta}:\RR\rightarrow\HH$ be a full connected preimage of $\beta$. Then there are lifts \pt of \p and
$\tilde{\delta}$ of $\delta$, so that $\tilde{\delta}$ has its initial point at \pt and its terminal point on $\tilde{\beta}$. Let 
$\tilde{\lambda}$ be the geodesic ray from \pt to $\tilde{\beta}_{+}$, one of the endpoints of $\tilde{\beta}.$ Then one  easily shows that for $g\in G$, either $g(\tilde{\lambda})\cap \tilde{\lambda}=\emptyset$ or else the intersection is the point set 
$\{\tilde{p}\}$ and $g$ is in the stabilizer of \pt. The ray $\tilde{\lambda}$ thus projects to a simple geodesic ray on \s with initial point \p. 

Observe that, given a sequence of values $t_j\rightarrow \infty,$ in the non-compact case $\lambda(t_j)$ has no limit points in $S$, while in the compact case the sequence $\lambda(t_j)$ will limit exactly at points on the geodesic $\beta.$

\subsubsection{Normalization of the covering group}
For the sake of clarity and computational simplicity we shall choose a nicely normalized Fuchsian group $G$ defining the orbifold \s.  
First, let the complex number $i\in \HH$ be our point \pt, covering the cone point \p. Then the transformation 
$$T_k (\zeta)= \frac{\cos\frac{\pi}{k}\zeta + \sin\frac{\pi}{k}}{- \sin\frac{\pi}{k}\zeta +\cos\frac{\pi}{k}}$$
generates the stabilizer of $i$ in $G$.

 We may further  normalize $G$ so that the geodesic $\lambda$ has a lift $\tilde{\lambda}_0$
with initial point i\ and terminal point $0\in \RR$. Then  $\tilde{\lambda}_1=T_k(\tilde{\lambda}_0)$ and 
 $\tilde{\lambda}_{-1}=T^{-1}_k(\tilde{\lambda}_0)$ are two other lifts of $\lambda$ beginning at $i$. The point $\tan\frac{\pi}{k}=a_k$ is the endpoint at infinity 
 of  $\tilde{\lambda}_1 $ and  $-a_k$ is the endpoint of $\tilde{\lambda}_{-1}.$ 
 
 We set some notation. Define the sets $I=[(-a_k,0)\times (0,\infty)]\cup [(0,a_k)\times (-\infty, 0)]$ and
 $J=[(-a_k,0)\times (a_k,\infty)]\cup [(0,a_k)\times (-\infty, -a_k)]$. We shall use the shorthand $\tilde{B}_r$ for the disc 
 $\tilde{B}_r(i)$ and $B_r$ for the disc $B_r(p).$ 
 
 \subsection{Excursions }
 \subsubsection{The definition}
 Given a cone neighborhood $B_r$ of \p and a geodesic ray \gam  on \s, we are concerned with  the "excursions" of \gam into $B_r$.   We saw in the introduction that these can be regarded as certain arcs of $\gamma,$ although 
 their definition is complicated by the self-overlap of $B_r$ for larger values of $r$. Hence, we would like to  distinguish the fine structure of the excursions as \gam intersects different overlapping pieces of $B_r$. 
 
 The above considerations dictate the 
 precise, but non-intrinsic, definition of an excursion which follows. Suppose $r>0.$ Define an {\it r-excursion} $e$ of a geodesic \gam into the cone neighborhood $B(r)$
to be   a lift \gamt of \gam to $\HH$ with 
 $( \tilde{\gamma}_{+}, \tilde{\gamma}_{-})\in I$ so that $\tilde{\gamma}\cap \tilde{B}_r\not =\emptyset.$

    The $r$-excursion $e=\tilde{\gamma}$ is called an {\it approximating } $r$-excursion  if $( \tilde{\gamma}_{+}, \tilde{\gamma}_{-})\in J.$ More generally, call $e$ an {\it approximating excursion} if $e$ is an approximating
 $r$-excursion for some $r\geq 0$. In Section \ref{section42}  we shall prove that every approximating excursion is an $r$-excursion for some $r< r_k= \sinh^{-1}(\cot\frac{\pi}{k}        ).$ 
  
 Let $E(\gamma, r)$ be the set of $r$-excursions of $\gamma$ and let $A(\gamma)$ be the set of approximating excursions of $\gamma$.

\subsubsection{Ordering the excursions along \gam }

To each $r$-excursion $e=\tilde{\gamma}$ there is a single intersection of the geodesics \gamt with the ray $\tilde{\lambda}_0$. We may therefore associate to $e$ the unique real parameter $t_e$, for which 
$\tilde{\gamma}(t_e)\in \tilde{\lambda}_0$.  Given a positive real value $r$ and a geodesic \gam the identification 
$e \rightarrow t_e$ defines a function  $\psi: E( \gamma, r)\rightarrow (-\infty ,\infty).$ Let $ E^{+}( \gamma, r)$ be the subset of $e\in  E( \gamma, r)$ so that $\psi(e)\geq 0$. $\psi$ also induces a map on $A(\gamma)$ and we similarly define  $A^{+}(\gamma).$ The range of $\psi$ shall be called the excursion parameters, written $\{t_e\}$. We then have

\begin{prop}
If $e$ and $e'$ are two distinct $r$-excursions along \gam, then $\psi(e)\not = \psi(e').$
\end{prop}

\begin{proof}
Write $e=\tilde{\gamma}$, $e= \tilde{\gamma}'$ and suppose $t_{e} = t_{e'}.$ Since $\tilde{\gamma}$ is distinct from $\tilde{\gamma}'$, there is a non-trivial
$g\in G$ so that $g(\tilde{\gamma})=\tilde{\gamma}'$ and $g(\tilde{\gamma}(t_{e}))=\tilde{\gamma}'(t_{e'}).$ It follows that $g(\tilde{\lambda}_0) \cap \tilde{\lambda}_0\not =\emptyset$. This is only possible if $g= T_k^n$ for some $n=1,\ldots, k-1.$ But then because we have $( \tilde{\gamma}_{+}, \tilde{\gamma}_{-})\in I,$ applying the transformation gives $(g( \tilde{\gamma}_{+}), g(\tilde{\gamma}_{-}))=( \tilde{\gamma}'_{+}, \tilde{\gamma}'_{-})\not\in I,$ which is impossible.
 
\end{proof}

Define an order on $E( \gamma, r)$ by stipulating that $e<e'$ if $\psi(e)< \psi(e').$ This ordering does not depend on the parametrization of $\gamma$. The next proposition shows that this ordering of $E$ is very well behaved.

\begin{prop}\label{prop2}
Given  a geodesic \gam, which is distinct from the simple closed geodesic $\beta$ used to define $\lambda$, suppose that the set of $r$-excursions $E^{+}=E^{+}( \gamma, r)$ is infinite.
Then there is a unique map from $\NN$ onto $E^{+}$, making $E^{+}$ into a sequence $\{e_j\}_{j=1}^{\infty}$ so that $e_j<e_k$ if and only if $j<k$. Furthermore, $\d{\lim_{j\rightarrow\infty}\psi(e_j)=\infty.}$
\end{prop}

\begin{proof}
As $G$ is discrete, $E^{+}$ must be countable. Therefore, it will suffice   to show that for any sequence $\{e_j\}_{j=1}^{\infty}$ in $E^{+},$ which is ordered as in the proposition, 
$ \lim_{j\rightarrow\infty}\psi(e_j)=\infty$. We argue by contradiction. Suppose there is such a sequence of $r$-excursions for which the excursion  parameters, which we shall write $\psi(e_j)=t_j$, do not diverge to $\infty.$ By passing to a subsequence, of the same name, we may suppose $t_j\rightarrow t^{*} <\infty.$ 

Write $e_j=\tilde{\gamma}_j.$ Choose a lift  $\tilde{\gamma}$   of $\gamma$. For each $j\in \NN$ there is a transformation $g_j\in G$ with $g_j(\tilde{\gamma}_j)=\tilde{\gamma}.$  These transformations must be distinct, since they respect the parameterization; in particular, $g_j(\tilde{\gamma}_j(0))=\tilde{\gamma}(0).$ Then $g_j(\tilde{B}_r)\cap \tilde{\gamma}\not =\emptyset$
  and $\tilde{\gamma}(t_j)=g_j(\tilde{\lambda}_0)\cap \tilde{\gamma}.$
We take the liberty to define, hopefully without confusion,  $g_j(\tilde{\lambda}_0)=\tilde{\lambda}_j.$ The rays  $\tilde{\lambda}_j$ intersecting  \gamt in a sequence of points $\tilde{\gamma}(t_j)$ converging to  $\tilde{\gamma}(t^{*}).$

 $\tilde{\lambda}_j$ has endpoints $g_j(i)\in \HH$ and $b_j\in \hat{\RR}.$ Since the interiors of the rays are all disjoint and the sequence   $\{t_j\}$ is increasing, the sequence $\{b_j\}$ can limit at, at most two points in $\in \hat{\RR},$ on opposite sides of \gamt. By restricting to one side, we may further stipulate that the sequence $\{e_j\}$ was chosen so that $b_j\rightarrow b^{*}\in \hat{\RR}.$ 

The initial point $g_j(i)$ of the ray  $\tilde{\lambda}_j$ is the center of the translate $g_j(\tilde{B}_r)$ of the disc $\tilde{B}_r.$
Since $G$ is discrete, the points $g_j(i)$ will limit on $\hat{\RR}$ and not in $\HH$. Consequently the discs $g_j(\tilde{B}_r)$ have Euclidean radii going to zero  and they also will limit on $\hat{\RR}$. Recall that for each $j$,  $g_j(\tilde{B}_r)\cap \tilde{\gamma}\not = \emptyset$. Therefore the points $g_j(i)$ all lie in the radius $r$ neighborhood of \gamt and must limit at an   endpoint of  $\tilde{\gamma},$ denoted by $\tilde{\gamma}_{*}$.

If the rays $\lambda_j$   accumulate in $\HH,$ it must be at a lift of $\beta.$ The limit point  $b_{*}$ cannot be an  endpoint of \gamt distinct from $\tilde{\gamma}_{*},$ for then the $\lambda_j$ would accumulate on \gamt, contradicting the hypothesis. Thus if $b^{*}$ is an endpoint of \gamt then $b_{*}=\tilde{\gamma}_{*}$.

Putting this all together: we have the sequence of geodesic rays $\tilde{\lambda}_j$, with the endpoints on $\hat{\RR}$ converging to $b^{*}$ and the endpoints in $\HH$ converging to $\tilde{\gamma}_{*}$. In order for this to happen the sequence of intersections $\tilde{\lambda}_j\cap\tilde{\gamma}=\tilde{\gamma}(t_j)$ must converge to $\tilde{\gamma}_{*},$
contrary to the assumption that the excursion parameters $\psi(e_j)=t_j$, do not diverge to $\infty.$ 

\end{proof}

  If the ray \gam contains infinitely many approximating excursions then they are ordered as a subsequence of the $r_k$-excursions. We write  $\{e_j^{*}\}_{j=1}^{\infty}$ for the sequence of approximating excursions in $A(\gamma)$.
 
 \subsubsection{The depth of an excursion}
 Given a geodesic $\alpha\in \HH,$ define $\tilde{d}(\alpha)$ to be the distance from $\alpha$ to the point $i$. 
 Define the   {\it depth} of the $r$-excursion $e=\tilde{\gamma}$ of \gam, to be the value $d(e)=\tilde{d}( \tilde{\gamma}).$   For a pair of points $(x,y)\in I,$ let $\alpha$ be the geodesic with endpoints $(\alpha_{+},\alpha_{-})=(x,y)$ and define 
 $D(x,y)= \tilde{d}(\alpha)$.

\section{The application of Ergodic Theory}\label{ergodic}

\subsection{The geodesic flow}
The unit tangent bundle $T_1\HH$ over $\HH$ may be identified with $\HH\times S^1.$ Define the measure $\tilde{\mu}=A\times\theta,$ where $A$ and $\theta$ are respectively, area measure on $\HH$ and lebesgue measure on $S^1.$ The measure $\tilde{\mu}$ is invariant under the geodesic flow $\tilde{G}^t$ on  $T_1\HH.
$ 

A full measure 
subset of $T_1\HH$ can be modeled by triples $(\psi, \zeta, t)\in \RR^3$ where $\psi \not = \zeta.$ Let $\alpha$ be a geodesic in $\HH$ parameterized so that  $\alpha(0)$ is the Euclidean midpoint of the semi-circle $\alpha(\RR)$. If $\psi=\alpha_{+}$ and $\zeta=\alpha_{-}$, then $(\psi, \zeta, t)$ corresponds to the unit tangent vector $\dot{\alpha}(t)\in T_1\HH.$ It is clear that the correspondence defines an injection onto the subset of non-vertical vectors in $T_1\HH$. In these coordinates the invariant measure  $\tilde{\mu}$ has the form $\d{\tilde{\mu}=\frac{1}{(\psi -\zeta)^2}d\psi d\zeta dt}$, up to a constant multiple.

The measure $\tilde{\mu}$ projects to a measure $\mu$ on $T_1S$, which is invariant under the the geodesic flow $G^t$ on $T_1S$. It is well known that when \s is a  finite area surface, $G^t$ acts ergodically with respect to $\mu$. Every orbifold has a finite branched cover that is
a hyperbolic surface. Thus ergodicity of the flow on the surface implies the ergodicity of the flow
on the orbifold. See   \cite{adler0},  \cite{nicholls} and \cite{series} for details.

\subsection{Convergence of the distribution}\label{section3}
Every $v_q \in T_1S$ uniquely determines a geodesic (or a geodesic ray) \gam where $\gamma(0)=q$ and  $\dot{\gamma}(0)=v_q,$ the unit vector tangent to \gam at $q$. Conversely, every geodesic determines a unit tangent vector at the parameter $t=0$. We say that a property holds for almost all geodesics on \s if there is a set $A\subset T_1S$, of full $\mu$-measure, so that the property holds for every geodesic determined by a vector in $A$, as above. 

There are two classes of subsets of $I$ and $J$ that play particularly important roles in what follows. For $0\leq z<\infty$ set
$\Omega(z)=\{(x,y)\in I\, |\, D(x,y)\leq z\}$ and for $0\leq z\leq r_k$ set $\Omega_{*}(z)=\{(x,y)\in J\, |\, D(x,y)\leq z\}.$ As a consequence of Corollary \ref{cor2}, $\Omega_*(z)=J$ for $z\geq r_k.$  At times we shall use the letters $\psi$ and $\zeta$ in place of $x$ and $y$. Then for appropriate values of $z$ define the integrals

\begin{equation*}\label{lambda1}
\Lambda(z)=\int_{\Omega(z)}\frac{1}{(\psi -\zeta)^2}d\psi d\zeta  \,\, \text{and}
\end{equation*}
\begin{equation*}\label{lambda2}
\Lambda^{*}(z)=\int_{\Omega_{*}(z)}\frac{1}{(\psi -\zeta)^2}d\psi d\zeta .
\end{equation*}

After the next theorem,   the proofs of  Theorems \ref{thm1} and \ref{thm2}
become computations in hyperbolic geometry  and plane integration.

\begin{thm}\label{reduction}
For almost all geodesic rays \gam on \s, for $r>0$ and   $0\leq z\leq r,\, dist_k(r,z)(\gamma)$ converges to $ dist_k(r,z),$ where 
\begin{equation*}
dist_k(r,z)=\frac{\Lambda(z)}{\Lambda(r)}.
\end{equation*}
For  almost all geodesic rays \gam on \s and for  $0\leq z\leq r_k,\, dist^{*}_k(z)(\gamma)$ converges to  $dist^{*}_k(z),$ where
\begin{equation*}
dist^{*}_k(z)=\frac{\Lambda^{*}(z)}{\Lambda^{*}(r_k)}.
\end{equation*}

\end{thm}

The remainder of   Section \ref{ergodic} shall be devoted to the proof of  this theorem. The idea is to produce a measure transverse to the geodesic flow, on which there is an ergodic action  related to the sequence of excursions along a generic geodesic on \s.  The distributions can be expressed as  limiting sums with respect to this action which  then, using the Birkhoff Ergodic Theorem, become the above integrals.

\subsection{A cross-section of $G^t$ over $T_1\lambda$}
Let  $T_1\lambda$ denote the subset of the unit tangent bundle over the geodesic $\lambda$, which consists of the vectors $v_q$ with $q\in \lambda\setminus \{p\}.$ We shall define  subsets of  $T_1\lambda$ associated with $r$-excursions and separately , with approximating excursions. Most of the work shall take place in $\HH,$ on the preimage $T_1\tilde{\lambda}_0$ of $T_1\lambda$. 

Given $r>0$ define the set of vectors $\tilde{l}^r\subset T_1\tilde{\lambda}_0$, where $v_q\in \tilde{l}^r$ if there is a geodesic $\alpha$ in $\HH$ with 

\begin{enumerate}
\item  $( \alpha _{+},  \alpha _{-})\in I,$ 
\item $\alpha(0)=q\in \tilde{\lambda}_0\setminus \{i\}$ and $\dot{\alpha}(0)=v_q,$
\item $\alpha\cap \tilde{B}_r\not =\emptyset,$ and
\item the geodesic ray $\pi\circ\alpha$ on \s contains infinitely many $r$-excursions.
\end{enumerate}

 Similarly, define $\tilde{l}^{*}\subset T_1\tilde{\lambda}_0$, where $v_q\in \tilde{l}^{*}$ if there is a geodesic $\alpha$ in $\HH$ with 
\begin{enumerate}
\item  $(  \alpha _{+},  \alpha _{-})\in J,$ 
\item $\alpha(0)=q\in \tilde{\lambda}_0 \setminus \{i\}$ and $\dot{\alpha}(0)=v_q,$ and
\item the geodesic ray $\pi\circ\alpha$ on \s contains infinitely many approximating excursions.
\end{enumerate}

In order to treat the various cases simultaneously, we adopt the convention that $\tau$ denotes either the parameter
$r>0$ or the symbol $*$. The set of vectors  $\tilde{l}^{\tau}$ projects to a set of vectors $l^{\tau}\subset T_1\lambda.$
Essentially, $ l ^{\tau}$ consists of the tangent vectors of the form $\dot{\gamma}(t_e),$ associated to  excursions $e$, with the caveat that the ray $\gamma$ contains infinitely many such excursions. Since \l is simple, the projection $\pi_{*}:T_1\tilde{\lambda}_0\rightarrow  T_1\lambda$ is bijective and consequently we have the following.

\begin{prop}
The projection $\pi_{*}:\tilde{l}^{\tau}\rightarrow l^{\tau}$ is a bijection.
\end{prop}

We let  $*$-excursion be another term for   approximating excursion.  
The next proposition asserts that  $l^{\tau}$ is a cross-section for the geodesic flow on \s. 

\begin{prop}\label{Mtau}
For each    value $\tau$, there is a full measure, $G^t$-invariant subset $M_{\tau}\subset T_1S$ so that 
$l^{\tau}$ is a cross-section for the flow $G^t$ acting on $M_{\tau}$. In particular, for each $v\in M_{\tau}$, there is a geodesic ray $\gamma$, possessing an infinite sequence of $\tau$-excursion parameters $\{t_{e_j}\}$, so that $\dot{\gamma}(0)=v$ and $G^{t}(v)\in l^{\tau}$ for $t>0$ if and only if $t=t_{e_j}$ for some $j\in \NN$.

\end{prop}
\begin{proof}
We argue when $\tau=r$. Let $M_r$ be the set of vectors $v\in T_1S$ so that a geodesic ray \gam with $\dot{\gamma}(0)=v$ contains infinitely many $r$-excursions, none of which intersect the cone point $p$. It follows from Proposition \ref{prop2} that $M_r$ is $G^t$-invariant. In particular, $l^r\subset M_r.$

In order to see that $M_r$ has full $\mu$-measure, choose a value $u$ so that $0<u<r_k$ and $u\leq r$. 
As a consequence of the Poincar\'{e} Recurrence Theorem \cite{sinai}, there is a set of full measure $N_u\subset T_1S$, with the property that for each $v\in N_u$ the geodesic \gam with $\dot{\gamma}(0)=v$ returns infinitely often to $B_u$.
It is proved in Corollary \ref{cor3} of Section \ref{section42}, that every intersection of a geodesic \gam with $B_u$ produces a $u$-excursion. Therefore, such a geodesic \gam must contain an infinite sequence of $u$-excursions. The subset $N_u^* $ consisting of $v\in N_u$ so that the $G^t$-orbit of $v$ does not pass thru $p$, has full measure in $N_u.$
 As $u\leq r$, $N_u^*\subset M_r$, showing that $M_r$ has full measure in $T_1S.$
 
Given $v\in M_r$, let \gam be the geodesic with $\dot{\gamma}(0)=v$. Then \gam contains an infinite sequence of $r$-excursions $\{e_j\}$ with excursion parameters $\psi(e_j)=t_{e_ j} .$ As before write $e_j= \tilde{\gamma}_j$ with $((\tilde{\gamma}_j)_{+},(\tilde{\gamma}_j)_{-})\in I.$ Thus for each $j\in \NN,\, 
G^{t_{e_j}}(v)=\dot{   \tilde{\gamma}}_j(t_{e_j})\in l^r.$

If for some $t>0$ and $v\in M_r\,,$ $G^t(v)\in l^r,$ then there is a lift $\tilde{v}$ of $v$ and \gamt of $\gamma$ so that  $\dot{\tilde{\gamma}}(t)=\tilde{v}\in T_1\tilde{\lambda}_0,$
$(\tilde{\gamma}_{+},\tilde{\gamma}_{-})\in I$ and  $ \tilde{\gamma}\cap \tilde{B}_r \not =\emptyset.$ Therefore $e=\tilde{\gamma}$ is an $r$-excursion.     The argument for approximating excursions is similar.  

\end{proof}

\subsection{The first return map}
Let $F_{\tau}:l^{\tau}\rightarrow l^{\tau}$ be the first return map under the geodesic flow $G^t : M_{\tau}\rightarrow M_{\tau}.$
In other words, given $v\in l^{\tau}$ let \gam be a geodesic in \s with $\dot{\gamma}(0) = v.$ Then $F_{\tau}(v)= \dot{\gamma}(t)$, where $t>0$ is the first value with $\dot{\gamma}(t)\in l^{\tau}.$ The invariant measure $\mu$ for the geodesic flow on \s induces an $F_{\tau}$-invariant measure $\nu_{\tau}$ on $l^{\tau}$ for which $F_{\tau}$ is ergodic, \cite{adler0,af0}.  Under the identification of $ l^{\tau}$  with $\tilde{l}^{\tau},$  $F_{\tau}$ lifts to  a map $\tilde{F}_{\tau}:\tilde{l}^{\tau}\rightarrow \tilde{l}^{\tau}$ which is invariant and ergodic with respect to the lifted measure $\tilde{\nu}_{\tau}$. Furthermore, this construction produces an isomorphism between the dynamical systems $(  l^{\tau}, F_{\tau}, \nu_{\tau})$ and $(\tilde{l}^{\tau}, \tilde{F}_{\tau}, \tilde{\nu}_{\tau})$. 

It will be easier to work with  $F_{\tau}$ in a form associated to the $(\psi,\zeta,t)$ coordinates.  Let $\Omega_{\tau}$ denote $\Omega(r)$ if $\tau=r$ and $J$ is $\tau=*$. There is a  map $B_{\tau}:\tilde{l}^{\tau}\rightarrow \Omega_{\tau},$
taking $v_q$ to the endpoints $(\alpha_{+}, \alpha_{-})$ of the geodesic $\alpha$ with  $\alpha(0)=q$ and
$\dot{\alpha}(0)=v_q$.  $B_{\tau}$ is injective in the complement of the subset of  $ \tilde{l}^{\tau}$  of periodic points of the map $\tilde{F}_{\tau}.$ Since the set of periodic points is countable, $B_{\tau}$ is injective on a set of full $\tilde{\nu}$-measure. Henceforth, we shall remove the set of periodic points from  $ \tilde{l}^{\tau}$ and 
$l^{\tau}.$

If we represent $v_q= (\psi_q,\zeta_q,t_q)$ in coordinates, then $B_{\tau}$ becomes the projection  $B_{\tau}(\psi_q,\zeta_q,t_q)=(\psi_q,\zeta_q),$ onto the first two coordinates. Let $\omega_{\tau}= B_{\tau}(\tilde{l}_{\tau}).$ and
define the transformations $R_{\tau}:\omega_{\tau}\rightarrow \omega_{\tau}$ by
 $R_{\tau}\circ B_{\tau}=B_{\tau}\circ \tilde{F}_{\tau}.$
Note that  $\omega_{\tau}$ is   precisely the subset of $\Omega_{\tau}$ so that if $\alpha$ is a geodesic with $(\alpha_{+}, \alpha_{-})=(\psi,\zeta)\in  \Omega_{\tau}$ the geodesic ray $\alpha$ projects to a geodesic ray on \s that contains infinitely many $r$-excursions or approximating excursions, not thru $p$, as appropriate.

Arguing as in \cite{series}, we can use a theorem of Ambrose \cite{amb} to show that the geodesic flow can be built up, in a nice way, from the first-return map on the cross section.  This theorem allows us to conclude that the $G^t$-invariant measure $\mu$ on \s has the form $\mu= \nu_{\tau}\times dt$ in local coordinates, where $dt$ represents lebesgue measure on the flow lines.
Then $\mu= \nu_{\tau}\times dt$ pulls back to the $\tilde{G}^t$-invariant measure $\tilde{\mu}=\tilde{\nu}_{\tau}\times dt$
on $T_1\HH.$  It follows that   in $(\psi,\zeta)$ coordinates, $\tilde{\nu}_{\tau}=\frac{1}{(\psi -\zeta)^2} d\psi d\zeta$ is the induced invariant measure for $R_{\tau}.$ We shall write $\tilde{\nu}_{\tau}=\tilde{\nu}$ since, in these coordinates, the measure is independent of $\tau.$

Since $M_{\tau} $ is full measure in $T_1S$ and $G^t$-invariant, the cross-section $l^{\tau}$ is of full measure in $T_1\lambda$. It follows, in particular, that $\omega_{\tau}$ is full measure in $\Omega_{\tau}.$ In light of this, the transformation $R_{\tau}:\Omega_{\tau}\rightarrow \Omega_{\tau}$ is well defined up to sets of measure zero, is $B_{\tau}$-conjugate to $\tilde{F}_{\tau},$ and has invariant measure $\tilde{\nu}_{\tau}.$ 

It is also important to note that, as a  consequence of Proposition \ref{Mtau}, for almost all $v\in T_1S$, there are infinitely many values $t_j\rightarrow \infty$ so that the geodesic \gam with $\dot{\gamma}(0)=v$ has $\dot{\gamma}(t_j)\in l^{\tau}$ . Therefore, for almost all $v\in T_1S$ the geodesic \gam with $\dot{\gamma}(0)=v$ has a lift with its endpoints in $\omega_{\tau}.$ 
In other words, almost every geodesic on \s has a lift with its endpoints in $\omega_{\tau}.$ 

\subsection{The proof of Theorem \ref{reduction}}
As usual, we shall prove the theorem for $r$-excursions. The details for approximating excursions are similar.  First observe that for $(\psi_0,\zeta_0)\in \Omega(r)$ we have the equality 
\begin{equation}\label{last}
\lim_{n\rightarrow\infty} \frac{1}{n}\#\{j\,|\,1\leq j\leq n , D\circ R ^j_r(\psi_0,\zeta_0)\leq z\}=
               \end{equation}       
   \begin{equation}\label{int}        
            \lim_{n\rightarrow\infty} \frac{1}{n}\sum_{j=1}^n(\chi_{[0,z]}\circ D )( R ^j_r(\psi_0,\zeta_0)).
     \end{equation}
The measure  $\tilde{\nu}^{*}=  \frac{1}{  \Lambda(r)  }  \tilde{\nu}  $ is an $ R $-invariant probability measure on $\Omega(r),$ with respect to which $\tilde{R}$ is ergodic.  By Birkhoff's Ergodic Theorem and the fact that $\omega_{r}$ has full $\tilde{\nu}^{*}$-measure in $\Omega(r)$, for almost all $(\psi_0,\zeta_0)\in \Omega(r)$ the limit (\ref{int}) is equal to 
                $$ \int_{\Omega(r)}(\chi_{[0,z]}\circ D)\tilde{\nu}^{*}                = \frac{1}{  \Lambda(r)  }\int_{\Omega(z)}\tilde{\nu}=\frac {\Lambda(z)}{  \Lambda(r)  }.$$
                
 Almost all geodesics \gam on \s have lifts with endpoints in $\omega_{\tau}$. Let \gam be one of these geodesics.   
 Then there is an infinite sequence     $e_j=\tilde{\gamma}_j,\, j=0,1,\ldots,$ of $r$-excursions along the ray \gam, with associated excursion parameters $t_j=\psi(e_j).$ Let $v_j=\dot{\gamma}(t_j)\in l^{\tau}.$ By Proposition \ref{Mtau}, $F_{\tau}(v_j)=v_{j+1}$ and therefore   $F^j_{\tau}(v_0)=v_j.$
The lift of $v_j$ to $\tilde{\lambda}_0$ is $\tilde{v}_j= \dot{\tilde{\gamma}}_j(t_j)$ and, as above, $\tilde{F}^j_{\tau}(\tilde{v}_0)=\tilde{v}_j.$

Write $\tilde{\gamma}=\tilde{\gamma}_0. $  To each $j \in \NN$ there is a unique $g_j\in G$ so that $\tilde{\gamma}_j=g_j(\tilde{\gamma}).$ Then in $(\psi,\zeta,t )$ coordinates, $\tilde{v}_j=((\tilde{\gamma}_j)_{+},(\tilde{\gamma}_j)_{-}, s_j)=(g_j(\tilde{\gamma}_{+}), g_j(\tilde{\gamma}_{-}),s_j)$, for some values $s_j$ that are not related to the $t_j$. Then if we set $(\psi_0,\zeta_0)=(\tilde{\gamma}_{+}, \tilde{\gamma}_{-})\in \omega_{\tau}$ we get
\begin{equation}\label{equalities}
 R_{\tau}^j(\psi_0,\zeta_0)= B_{\tau}\circ\tilde{F}_{\tau}^j\circ B_{\tau}^{-1}(\tilde{\gamma}_{+}, \tilde{\gamma}_{-})=
B_{\tau}\circ\tilde{F}_{\tau}^j(\tilde{v}_0)=
\end{equation}
$$
 B_{\tau}(\tilde{v}_j)= B_{\tau}(g_j(\psi_0), g_j(\zeta_0),s_j)=(g_j(\psi_0), g_j(\zeta_0)). $$

Let $\chi_Y$ denote the characteristic function of the set $Y$. Turning to the distribution we have
$$dist_k(r,z)(\gamma)= \lim_{n\rightarrow\infty} \frac{1}{n}\#\{j\,|\,1\leq  j\leq n , d(e_j)\leq z\}= $$

                  $$ \lim_{n\rightarrow\infty} \frac{1}{n}\#\{j\,|\,1\leq j\leq n , \tilde{d}(\tilde{\gamma}_j)\leq z\}$$
                  
           $$  \lim_{n\rightarrow\infty} \frac{1}{n}\#\{j\,|\,1\leq j\leq n , D(g_j(\tilde{\gamma}_{+}), g_j(\tilde{\gamma}_{-}))\leq z\}= $$  
   \begin{equation}\label{now}
\lim_{n\rightarrow\infty} \frac{1}{n}\#\{j\,|\,1\leq j\leq n , D(g_j(\psi_0),g_j(\zeta_0))\leq z\}.
\end{equation}
By the sequence of equalities (\ref{equalities}), the limit (\ref{now}) is equal to the limit (\ref{last}). Therefore, we may conclude that for almost all geodesics \gam on \s, $dist(r,z)(\gamma)$ converges to 
$dist(r,z)=\frac{\Lambda(z)}{\Lambda(r)},$ which proves Theorem \ref{reduction} for $r$-excursions.
$ \hfill\Box  $

\section{The computation of $\Lambda$ and $\Lambda^{*}$} \label{section40}
\subsection{The values of  $\Lambda$ and $\Lambda^{*}$}

In light of Theorem \ref{reduction}, the proofs of Theorems \ref{thm1} and \ref{thm2} will be completed by computing the values of  $\Lambda$ and $\Lambda^{*}$.     
Recall that $r_k= \sinh^{-1}(\cot\frac{\pi}{k}        )$, $\delta_k=\sinh^{-1}(\cot  \frac{2\pi}{k}  ) $ and 
$\varphi_k(x)= \frac{1-x \tan  \frac{\pi}{k}}{x+ \tan  \frac{\pi}{k}}.$
   
     \begin{thm}\label{thmlambda}
For $z\geq 0$ we have

$$\Lambda(z)=
\left\{\begin{array}{cc}
     \frac{2\pi}{k}\sinh z  
 & {\rm if} \, z\leq\, r_k\\
 \\
   \sinh z \tan^{-1}(\frac{1}{\sinh z}) + \log(\sin  \frac{\pi}{k}\cosh z) 
  & {\rm if} \, z> r_k , 
     \end{array} \right.
     $$
$$\Lambda^{*}(z)=
\left\{\begin{array}{cc}
     \frac{2\pi}{k}\sinh z  
 & {\rm if} \, z\leq\,\delta_k\\
 \\
 \sinh z\tan^{-1}(\varphi_k(\sinh z)) +\log(2\cosh z \sin \frac{\pi}{k}\cos\frac{\pi}{k})   
 & {\rm if} \, \delta_k< z< r_k\\
 \\
 \log (2 \cos  \frac{\pi}{k}) 
 & {\rm if} \, z\geq\, r_k.
\end{array} \right.
     $$
When $k=3$ or 4, $\delta_k\leq 0$ and only the second two cases occur for the value of $\Lambda^{*}(z)$.
\end{thm}

The rest of this section is devoted to the proof of this theorem.

\subsection{Hyperbolic geometric considerations} 
 \subsubsection{The geodesic tangent to a disc}\label{section4}
In order to turn $\Lambda$ and $\Lambda^{*}$ into  easily computable  double integrals, we shall derive several formulae that will be of use in determining the limits of integration.  Given distinct points $a,b\in \hat{\RR}$, let $\overline{ab}$ denote the geodesic $\alpha$ in $\HH$ with $\alpha_{+}=a$ and $\alpha_-=b$.

\begin{thm}\label{wrho}
Given $x\geq0$, the point $w\in [-\frac{1}{x}, x)$ for which the geodesic $\overline{xw}$, is tangent to the hyperbolic disc $\tilde{B}_{\rho} $ is given by the formula 
$$w=W_{\rho}(x)=\d{\frac{x\sinh\rho -1}{x+\sinh\rho}}.$$
Similarly, given $x\leq0$ the point $w\in (x,-\frac{1}{x}]$ for which the geodesic $\overline{xw}$, is tangent to the hyperbolic disc $\tilde{B}_{\rho}$ is given by the formula 
$w=-W_{\rho}(-x) .$

\end{thm}

Let $\Delta$ denote the unit disc model for the hyperbolic plane, which we envision as sitting in the complex plane.  Write $\partial \Delta$ for the unit circle, which is its boundary at infinity. The proof of Theorem \ref{wrho} involves computations in $\Delta$ and the transformation of those results to $\HH$. The following lemma can be proved by an easy computation using the hyperbolic metric in $\Delta.$ 

\begin{lemma}\label{lem1}
In $\Delta,$ the circle of hyperbolic radius $\rho$ and center 0 has Euclidean radius
$r= \tanh\frac{\rho}{2}$
\end{lemma}

If $z, \xi \in \partial \Delta$  then, as in the previous case, we let 
$\overline{z\xi}$ be the hyperbolic geodesic in $\Delta$ with endpoints $z, \xi$. Similarly, for  $z\in \partial \Delta$ let $\overline{z0}$ be the geodesic with endpoint $z$ passing thru the origin. This geodesic is a Euclidean straight line. We use the same notation $\tilde{B}_{\rho}(0)$ to denote the hyperbolic disc of center 0 and radius $\rho$ in $\Delta.$ To avoid confusion we shall continue to write   $\tilde{B}_{\rho}$ for    $\tilde{B}_{\rho}(i)$.

\begin{lemma}\label{lemma2}
Suppose that $z, \xi\in\partial \Delta$ are two points so that the geodesic $\overline{z\xi}$ is tangent to the the hyperbolic disc
$\tilde{B}_{\rho}(0)$. Then $|z-\xi|=2\text{sech}\rho.$
\end{lemma}

\begin{proof}
There is no loss of generality in supposing that $z=x+iy$ with $x,y\geq 0$ and $\xi = \overline{z}$.
Then we have $|z-\xi|=2y$ and it remains to show that $y= \text{sech}\rho.$
Let $C$ be the circle perpendicular to $\partial \Delta$ whose intersection with $\Delta$ is the geodesic $\overline{z\xi}$.
Since $C$ is invariant under the reflection $z\rightarrow \frac{1}{\overline{z}}$ which fixes $\partial \Delta$, it crosses the real axis in points $0<r<1 $ and $\frac{1}{r}.$ By Lemma \ref{lem1}, since $C$ is tangent to the circle of radius $\rho,\, r=\tanh\frac{\rho}{2}.$

As a consequence of $z$ lying on $C,$ it satisfies the equation $|z-\frac{1}{2}(\frac{1}{r}+r)|^2=[\frac{1}{2}(\frac{1}{r}-r)]^2.$
Also, since $z \in\partial \Delta$, $|z|=1$. Solving these equations simultaneously gives $y=\frac{1-r^2}{1+r^2}.$
 Setting  $ r=\tanh\frac{\rho}{2}$ and simplifying, gives the result.
\end{proof}

\begin{lemma}\label{previous}
Suppose $z$ is a point in the upper half of the unit circle $\partial \Delta$. Let $\xi\in \partial \Delta$ be the point which lies to the right of the line $\overline{z0}$ so that $\overline{z\xi}$ is tangent to the disc $\tilde{B}_{\rho}(0)$. Then
$\xi = z(c-i\sqrt{1-c^2})$ where $c= 1-2\text{sech}^2\rho.$ The above also holds for $z=1,$ where $\xi$ is then a point in the lower half of $\partial \Delta.$
\end{lemma}

\begin{proof}
From the previous lemma we have $|z-\xi|=2\text{sech}\rho,$ which simplifies to 
$\text{Re} z\overline{\xi}= 1-2\text{sech}^2\rho.$ For simplicity write $c=1-2\text{sech}^2\rho.$ Note that as $\rho$ increases on the interval $[0,\infty),\,$ $c$ increases on $[-1,1)$. 

The point $\xi$ may be written in the form $\xi= \eta z$, for a  point $\eta$ of modulus one.
Let $\eta=u+iv$. One see by a simple computation that $u=\text{Re}z\overline{\xi}=c$ and then $w=\pm\sqrt{1-c^2}$.
With the choice of the minus sign the points $\eta = c-i\sqrt{1-c^2}$ fill out the lower half of the unit circle. As a function of $\rho,\,$ $\eta$ moves counter-clockwise through its arc as $\rho$ increases from zero to infinity. Thus, treating the point $\xi=\eta z$ as a function of $\rho$, we see that it moves counter-clockwise around the circle from $-z$ to $z$ as   $\rho$ increases from zero to infinity.  It follows that for all values of $\rho$, the point $\xi=\eta z$ lies to the right of the line $\overline{z0},$ or in the lower half-plane if $z=1.$ With the choice of a plus sign, $\eta$ fills out the upper half of the circle and the points $\eta z$ lie to the left of $\overline{z0},$ or in the upper half-plane if $z=1.$  That completes the proof.
\end{proof}

\noindent {\em Proof of Theorem \ref{wrho}.}\,
The M\"obius transformation $g(z)=\frac{-iz+i}{z+1}$ maps $\Delta$ isometrically to $\HH,$ taking the disc $\tilde{B}_{\rho}(0)$ to the disc $\tilde{B}_{\rho} $. $g$ also defines a 
 one-to-one correspondence between $\partial \Delta$ and $ \hat{\RR}$, where 1 goes to 0, -1 goes to $\infty$ and the upper and lower semi-circles   respectively, are mapped to the positive and negative real axes. 
 
 As in the previous argument let  $\eta = c-i\sqrt{1-c^2},$ where $c= 1-2\text{sech}^2\rho,$ and define the transformation $W_{\rho}(x)=g(\eta g^{-1}(x))$.  Simplifying we get 
 $$W_{\rho}(x)=g(\eta g^{-1}(x))=
 \frac{ix\left(\frac{1+\eta}{1-\eta}\right)-1}{x+i \left(\frac{1+\eta}{1-\eta}\right) }.$$
 A further computation gives 
 $$i\left(\frac{1+\eta}{1-\eta}\right)=      -\frac{2\text{Im}\eta}{|1-\eta |^2}=
 \frac{\sqrt{1-c^2}}{1-c}=\sinh\rho,$$
  which shows that the  formula for $W_{\rho}$ is the one asserted in the theorem.
 Now we need to see that $W_{\rho}$ does what it is claimed to do.

 If $x\geq 0$ then  $g^{-1}(x)=z$ lies in the upper-half of $\partial \Delta$ or is equal to 1. Then $\xi=\eta z$ is the point of $\partial \Delta$ for which  the geodesic $\overline {z\xi}$ satisfies the hypotheses of the previous lemma.  The image of this geodesic, $g(\overline {z\xi}),$ is the hyperbolic geodesic in $\HH$ which is tangent to $g(\tilde{B}_{\rho}(0))=\tilde{B}_{\rho} $ and has  endpoints $x$ and $W_{\rho}(x)$. Recall from the proof of Lemma \ref{previous} that $\xi=\eta z$ lies in the   counter-clockwise   arc of $\partial \Delta$  between $-z$ and $z$. Since the image under $g$ of this arc is the interval $[-\frac{1}{x},x)$, the proof is complete for $x\geq 0$.
 
Now suppose $x\leq 0.$ Consider the isometry $h(z)=-\overline{z}$ of $\HH.$ $h$ maps the geodesic tangent to $\tilde{B}_{\rho} $, with endpoints $-x$ and $W_{\rho}(-x)\in [\frac{1}{x} ,-x)$ to the geodesic tangent to $\tilde{B}_{\rho} $ with endpoints $x$ and  $-W_{\rho}(-x)\in (x,-\frac{1}{x} ]. $ Thus  $-W_{\rho}(-x)$ is as asserted in the proposition.
$ \hfill\Box $ 

\subsubsection{Consequences of the tangency computations} \label{section42}
 \begin{cor}\label{cor1}
 Given $x,z\geq 0,\, D(x,y)\leq z$ for $y\in [-\frac{1}{x},x)$ if and only if   $y\in [-\frac{1}{x},W_z(x)].$ 
 \end{cor}
 
 \begin{proof}
 For fixed $x\geq 0$ and for $y\in [-\frac{1}{x},x)\,$,  define $  H(y) =D(x,y).$  Then $  H(y)$ is the value $\rho$ for which the geodesic $\overline{xy}$ is tangent to the disc 
 $\tilde{B}_{\rho} $. We have shown in the proofs of Lemma \ref{previous} and Theorem \ref{wrho} that $H$ is an increasing function of $y$, mapping $ [-\frac{1}{x},x)$ onto $[0,\infty)$.
 Also recall that, for $x,y$ as above, $\overline{xy}$ is tangent to  $\tilde{B}_{z} $ precisely when $y=W_z(x).$  
It follows that $D(x,y)=H(y)\leq z$ if and only if  $y\in [-\frac{1}{x},W_z(x)].$
 \end{proof}
 
 \begin{cor}\label{cor2}
 When $k\geq 5$, every $z$-excursion with $z\leq\delta_k$ is an approximating excursion. Moreover, every approximating excursion is a
 $z$-excursion for some $z<r_a$. 
 \end{cor}

 \begin{proof}
 First observe that the equation $W_z(a_k)=-a_k$
  has the solution 
 $$
  \begin{array}{cc}
  z=\sinh^{-1}(\varphi_k^{-1}(a_k))
  =\sinh^{-1}(\frac{1-\tan^2\frac{\pi}{k}}{2\tan  \frac{\pi}{k}})=
 \sinh^{-1}(\cot \frac{2\pi}{k})=\delta_k.
 \end{array}
 $$
 This is positive if and only if $k\geq 5,$ which we henceforth suppose to be the case. 
  For $x,z\geq 0,\, W_z(x)$ is an increasing function in both variables. Therefore for $z\leq \delta_k$ and $x\in (0,a_k),$ $W_z(x)<W_z(a_k)\leq W_{\delta_k}(a_k)=-a_k.$
 
Let $e=\tilde{\gamma}$ be a $z$-excursion for some $z\leq \delta_k.$ Then $x=  \tilde{\gamma}_+\in (0,a_k)$ and, as a consequence of the above, $ \tilde{\gamma}_-=W_z(x)<-a_k.$ This says that $e$ is an approximating excursion, proving the first assertion of the corollary.
 
Now observe that $z=D(0,-a_k)$ is the value for which $W_z(0)=-a_k.$  Solving for $z$ we get 
$z=\sinh(\cot(\frac{\pi}{k}))=r_k.$ Thus the geodesic $\alpha$ with endpoints 0 and $a_k$ is tangent to the disc $\tilde{B}_{r_k}.$

Let $e=\tilde{\gamma} $ be an approximating excursion. Then $ \tilde{\gamma}_+\in (0,a_k)$ and 
$ \tilde{\gamma}_ -\in [-\frac{1}{x},-a_k).$ It follows that \gamt lies above $\alpha$ in $\HH$ and therefore,
$\tilde{\gamma}\cap \tilde{B}_{r_k}\not = \emptyset.$  The second part of the corollary is proved.

\end{proof}

\begin{cor}\label{cor3}
Suppose \gam is a geodesic on \s and $\gamma\cap B_{r_k}=\epsilon\not = \emptyset.$ 
Then there is  a lift \gamt of \gam and a lift $\tilde{\epsilon}$ of $\epsilon$ contained in \gamt, so that $\tilde{\epsilon}\cap \tilde{B}_{r_k}\not = \emptyset$  and
$e=\tilde{\gamma}$ is a
$z$-excursion for some $z<r_k.$
\end{cor}

\begin{proof}
 
Choose a lift \gamt so that $\tilde{\epsilon}\cap \tilde{B}_{r_k}\not = \emptyset$  and $\tilde{\gamma}_+\in (0,a_k).$ Without loss of generality we suppose that $\tilde{\gamma}_-< \tilde{\gamma}_+. $ Since $W_{r_k}(a_k)=0,$ the geodesic $\alpha$ with endpoints 0 and $a_k$ is tangent to the disc $ \tilde{B}_{r_k}$. If $\tilde{\gamma}_->0$ then \gamt lies entirely below $\alpha$ and $\tilde{\gamma}\cap \tilde{B}_{r_k}\not = \emptyset,$ which is impossible. Therefore  $\tilde{\gamma}_-<0$ and $e=\tilde{\gamma}$ is a
$z$-excursion for some $z<r_k.$
\end{proof}

 \subsection{The proof of Theorem \ref{thmlambda}}
 We shall write the integrals as sums of iterated double integrals. Once in this form, their actual  evaluation is a elementary.
 
\subsubsection{The computation of $\Lambda(z)$}

 Consider fixed values $z\geq 0$ and $x\in (0,a_k)$. If $W_z(x)<0$ then, by Corollary \ref{cor1},  $D(x,y)\leq z$ for $y\in [-\frac{1}{x},0)$ precisely when 
 $y\in [ -\frac{1}{x}, W_z(x)]$, whereas if  $W_z(x)\geq 0$ then $D(x,y)\leq z$ for all   $y\in [ -\frac{1}{x}, 0).$ As a result of this observation, for fixed $z>0$, there are two cases to be considered.

The first case is when $z\leq r_k$. This condition is equivalent to   $a_k \leq \frac{1}{\sinh z}$  or  $W_z(a_k)\leq 0.$ Since $W_z(x)$ is an increasing function of $x$, this occurs if and only if $W_z(x)<0$ for all $x\in (0,a_k).$ Then, as observed above,
  for  $x\in (0,a_k),\, D(x,y)\leq z$ for $y\in [ -\frac{1}{x},0)$   if and only if $y\in [ -\frac{1}{x}, W_z(x)].$  By symmetry, we have that for  $x\in (-a_k,0),\,  D(x,y)\leq z$ for $y\in (0, -\frac{1}{x}]$ if and only if $y\in [-W_z(-x),-\frac{1}{x}] .$
 Consequently,  when $z\leq r_k$ we get
  \begin{equation*}
   \d{ \Lambda(z)=\int_0^{a_k}  \int_{-\frac{1}{x}}^{W_z(x)}\frac{1}{(x -y)^2} dy dx  \,+\,
 \int_0^{-a_k}  \int^{-\frac{1}{x}}_{-W_z(-x)}\frac{1}{(x -y)^2} dy dx } 
\end{equation*}

In the second case $z>r_k$ and we will have $W_z(x)= 0$ for   $x= \frac{1}{\sinh z}\in (0,a_k).$ Then  the set of $y$ values for which  $D(x,y)\leq z$ will take two different forms. For $x\in (0,  \frac{1}{\sinh z}),$ we have $ D(x,y) \leq z$ for $y\in [ -\frac{1}{x},0)$  if and only if  $y\in [ -\frac{1}{x}, W_z(x)];$ whereas,  for $x\in [ \frac{1}{\sinh z},0)$ we get $ D(x,y) \leq z$ for all $y\in [-\frac{1}{x},0).$ Similar results again hold, by symmetry, for $x\in (-a,0).$ 
 Then for  $z>r_k$ we get
  \begin{equation*}
 \begin{split}
  \Lambda(z)=&\d{\int_0^{\frac{1}{\sinh z}}  \int_{-\frac{1}{x}}^{W_z(x)}\frac{1}{(x -y)^2} dy dx  \,+\, 
 \int^{a_k}_{\frac{1}{\sinh z}}  \int_{-\frac{1}{x}}^{0}\frac{1}{(x -y)^2} dy dx   \,+}\\
 & \d{ \int^0_{-\frac{1}{\sinh z}}  \int^{-\frac{1}{x}}_{-W_z(-x)}\frac{1}{(x -y)^2} dy dx  \,+\, 
 \int_{-a_k}^{-\frac{1}{\sinh z}}  \int^{-\frac{1}{x}}_{0}\frac{1}{(x -y)^2} dy dx  }\\
  \end{split}
\end{equation*}
 \subsubsection{The computation of $\Lambda^{*}(z)$}
 The same considerations and cases come into play in the computation of $\Lambda^{*}(z)$. In  addition there is a third special case,   which is a consequence of the fact that every approximating excursion is an  $r$-excursion for   some $r<r_k.$  The three cases depend on the position of the points $W_z(x)$ with respect to $-a_k$, rather than with respect to 0, as  in the earlier computations. The cases are also described by the position of the particular point $x_z$, for which  $W_z(x_z)=-a_k.$  
 
  In the first case   $ z\leq \delta_k$. As observed in the proof of Corollary \ref{cor2}, this is  equivalent to $W_z(x)<-a_k$ for all $x\in (0,a_k)$.  Since $\delta_k\leq 0$ for $k=3,4$, this case does not arise for those values of $k$. Then, by Corollary \ref{cor1},  for all $x\in (0,a_k),\, D(x,y)\leq z$ for $y\in [ -\frac{1}{x},-a_k)$ if and only if $y\in [ -\frac{1}{x}, W_z(x)],$  and by symmetry, we have that for all $x\in (-a_k,0),\,  D(x,y)\leq z$ for $y\in ( a_k, -\frac{1}{x}]$ if and only if $y\in [-W_z(-x),-\frac{1}{x}] .$ This is the same situation we had computing the integral in the earlier first case. Therefore for $z\leq \delta_k,$  
  $$ \Lambda(z)=\Lambda^{*}(z)=\frac {2\pi}{k}\sinh z.$$  
 
 In the second case   $\delta_k<z<r_k,$ where we let $\delta_k=0$ for $k=3,4$.  Consider the equation 
 $W_z(x)=-a_k,$ whose solution is a point $x_z=\varphi_k(\sinh z).$ $x_z$ is a decreasing function of $z$. We've seen that $x_{\delta_k}=a_k$ and $x_{r_k}=0$. Therefore, for $\delta_k<z<r_k$ we have $W_z(x)=-a_k$ for the point $x=x_z\in (0,a_k).$  
    As  in case two in the computation of $\Lambda$, the set of $y$ values for which  $D(x,y)\leq z$ will take two different forms. For $x\in (0, x_z),$ we have $ D(x,y) \leq z$ for $y\in [-\frac{1}{x},-a_k)$ if and only if  $y\in [ -\frac{1}{x}, W_z(x)]$ and for $x\in [ x_z,0),$ we get $ D(x,y) \leq z$ for all $y\in [-\frac{1}{x},-a_k).$ Similar results again hold, by symmetry, for $x\in (-a_k,0).$ 
 It follows that for $\delta_k<z<r_k$
\begin{equation*}
 \begin{split}
  \Lambda^{*}(z)=&\d{\int_0^{x_z}  \int_{-\frac{1}{x}}^{W_z(x)}\frac{1}{(x -y)^2} dy dx  \,+\, 
 \int^{a_{k}}_{x_z}  \int_{-\frac{1}{x}}^{0}\frac{1}{(x -y)^2} dy dx   \,+}\\
 &\d{ \int^0_{-x_z}  \int^{-\frac{1}{x}}_{-W_z(-x)}\frac{1}{(x -y)^2} dy dx  \,+\, 
 \int_{-a_k}^{-x_z}  \int^{-\frac{1}{x}}_{0}\frac{1}{(x -y)^2} dy dx  }\\
  \end{split}
\end{equation*}

The last case to consider is when $z\geq  r_k$. By Corollary \ref{cor2}, $\Omega_*(z)= J$
for all $z\geq  r_k.$  Therefore, for all $x\in (0,a_k),\, D(x,y)\leq z$ for all $y\in [ -\frac{1}{x},-a_k).$ Consequently, in this last case we have, for $z\geq r_k$
 \begin{equation*} 
\begin{split}
& \d{\Lambda^{*}(z)= \int_0^{a_k}  \int_{-\frac{1}{x}}^{-a_k}\frac{1}{(x -y)^2} dy dx  \,+\,
 \int_0^{-a_k}  \int^{-\frac{1}{x}}_{a_k}\frac{1}{(x -y)^2} dy dx}
 \end{split}
 \end{equation*}
 
\section{Depth as a function of area}\label{section5}

We have been dealing with  excursion depth as expressed in  terms of the radius of the largest cone neighborhood disjoint from the excursion. In \cite{haas}, where we   looked at excursion  into a cusp, it was just not possible to express the depth of an excursion as the distance to the cusp end, since the cusp is off at infinity. Instead we considered the area of the largest neighborhood of the cusp disjoint from the excursion. Suppose we were to take the same approach here. 

In the hyperbolic plane $\HH$ the area of a hyperbolic disc of radius $r$ is $2\pi (\cosh r -1)$ \cite{Beardon}. Let $G$ be a Fuchsian group representing \s and let $\pi :\HH \rightarrow \HH/G=S$ denote the natural projection map. Suppose \p is the cone point on \s and \pt is a preimage in $\HH$. The actual area of a disc $B_r(p)$   is at least $1/k$ times 
  the area of its preimage $\tilde{B}_r(\tilde{p})$, with equality if the projection is precisely k-to-1. In the case of equality, if $A$ is the area of $B_r(p)$, then
   $\tilde{B}_r(\tilde{p})$ has area $kA= 2\pi (\cosh r -1)$.  Even if  there is overlap and the projection is not k-to-1, in order to be consistent with \cite{haas}, we define the A-depth of an excursion $e$, written $A(e)$, to be $\frac{2\pi}{k} (\cosh (d(e))-1)$.
  Define   an $R$-excursion $e$ to be a usual $r$-excursion where $r=\cosh^{-1}(\frac{k}{2\pi}R+1)$, the radius of the disc of area $R$ about \p. The definition of an approximating excursion is unchanged. 
 
 Suppose $\{e_j\}$ is a sequence of $R$-excursions along a geodesic $\gamma$ on \s. Then define
 $$ Adist_k(R,Z)(\gamma)= \lim_{n\rightarrow\infty} \frac{1}{n}\#\{j\,|\,1\leq j\leq n , A(e_j)\leq Z\}.$$
 As a consequence of the definition of $A(e)$, this is equal to 
 $$ \lim_{n\rightarrow\infty} \frac{1}{n}\#\{j\,|\, 1\leq j\leq n , d(e_j)\leq\cosh^{-1}(\frac{k}{2\pi}Z+1)\}.
  $$
  For almost all geodesics \gam on \s this last limit was shown to converge to
$$ dist_k(\cosh^{-1}(\frac{k}{2\pi}R+1),\cosh^{-1}(\frac{k}{2\pi}Z+1))=
\frac{\Lambda(\cosh^{-1}(\frac{k}{2\pi}Z+1))}{\Lambda( \cosh^{-1}(\frac{k}{2\pi}R+1))}.$$
 
For $R \leq \frac{2\pi}{k} (\cosh (r_k)-1) $ and $Z\leq R$, the above takes the form
\begin{equation*}
\sqrt{\frac{ (Z+\frac{2\pi}{k})^2-(\frac{2\pi}{k})^2}{ (R+\frac{2\pi}{k})^2-(\frac{2\pi}{k})^2}}, 
\end{equation*}
which is almost linear in $Z$ and which approaches $Z/R$ as $k\rightarrow \infty.$ This last expression in the variable $Z$ is the  
  distribution function for an $R$-excursions along a geodesic with respect to  a cusped end of a surface, when $R\leq 2,$ \cite{haas}.
Note also that as $k \rightarrow\infty,\,$ $R_k$ approaches 2. 

For approximating excursion, the corresponding  distribution in terms of area has the form
$$\frac{\sqrt{ (Z+\frac{2\pi}{k})^2-(\frac{2\pi}{k})^2}}{2\log(2\cos\frac{\pi}{k})},$$
for $Z\leq \frac{2\pi}{k} (\cosh (\delta_k)-1) $. Again as $k\rightarrow\infty$ the distributions limit at 
$Z/2\log 2$, the value of the distribution of approximating excursions alone geodesics out a cusped end when 
$Z\leq 1.$

Clearly, the other cases look rather  nasty. Nevertheless, one can verify that they limit at the corresponding distributions
on surfaces with cusped ends.


\end{document}